\input amstex
\documentstyle{amsppt}
%\magnification 1200
\NoBlackBoxes

\TagsOnRight

\def\cal{\Cal}
\def\AA{{\cal A}}

\def\HH{{\cal H}}
\def\MM{{\cal M}}

\def\LL{{\cal L}}

\def\Z{{\Bbb Z}}
\def\C{{\Bbb C}}
\def\R{{\Bbb R}}

\def\Q{{\Bbb Q}}
\def\e{{\epsilon}}

\def\n{\noindent}
\def\part{{\partial}}
\def\dudtau{{\part u \over \part \tau}}

\rightheadtext{Hamiltonian diffeomorphism group} \leftheadtext{
Yong-Geun Oh }

\topmatter
\title
Spectral invariants and length minimizing property of Hamiltonian
paths
\endtitle
\author
Yong-Geun Oh\footnote{Partially supported by the NSF Grant \#
DMS-0203593,  Vilas Research Award of University of Wisconsin and
by a grant of the Korean Young Scientist Prize \hskip8.5cm\hfill}
\endauthor
\address
Department of Mathematics, University of Wisconsin, Madison, WI
53706, ~USA \& Korea Institute for Advanced Study, 207-43
Cheongryangri-dong Dongdaemun-gu Seoul 130-012, KOREA
\endaddress

\abstract In this paper we provide a criterion for the
quasi-autonomous Hamiltonian path (``Hofer's geodesic'') on
arbitrary closed symplectic manifolds $(M,\omega)$ to be length
minimizing in its homotopy class in terms of the spectral
invariants $\rho(G;1)$ that the author has recently constructed.
As an application, we prove that any
autonomous Hamiltonian path on {\it arbitrary}
closed symplectic manifolds is length minimizing in
{\it its homotopy class} with fixed ends, as long as
it has no contractible periodic orbits {\it of period one}
and it has a maximum and a minimum that are {\it generically
under-twisted}, and all of its critical points are non-degenerate
in the Floer theoretic sense.

\endabstract


\endtopmatter

\document

\bigskip

\centerline{\bf Contents} \medskip

\n \S1. Introduction \smallskip

\n \S2. Preliminaries \smallskip

\n \S3. Chain level Floer theory and spectral invariants
\smallskip

\n \S4. Construction of fundamental Floer Novikov cycles
\smallskip

\n \S5. The case of autonomous Hamiltonians
\smallskip

\smallskip

\bigskip

\head  \bf \S1. Introduction
\endhead
In [H1], Hofer introduced an invariant pseudo-norm on the group
$\HH am(M,\omega)$ of compactly supported Hamiltonian
diffeomorphisms of the symplectic manifold $(M,\omega)$ by putting
$$
\|\phi\| = \inf_{H\mapsto \phi} \|H\| \tag 1.1
$$
where $H \mapsto \phi$ means that $\phi= \phi_H^1$ is the time-one
map of Hamilton's equation
$$
\dot x = X_H(x),
$$
and $\|H\|$ is the function defined by
$$
\|H\| = \int _0^1 \text{osc }H_t \, dt = \int_0^1(\max H_t -\min
H_t)\, dt
$$
which is the Finsler length of the path $t \mapsto \phi_H^t$. He
[H2] also proved that the path of any compactly supported {\it
autonomous} Hamiltonian on $\C^n$ is length minimizing as long as
the corresponding Hamilton's equation has no non-constant periodic
orbit {\it of period less than or equal to one}. This result has
been generalized in [En], [MS] and [Oh3] under the additional
hypothesis that the linearized
flow at each fixed point is not over-twisted i.e., has no closed
trajectory of period less than one. The latter hypothesis was
shown to be necessary for any length minimizing geodesics with
some regularity condition on the Hamiltonian path [U], [LM].
The following result is the main result from
[Oh3] restricted to the autonomous Hamiltonians among other results.

\proclaim{Theorem I [Oh3]}
Let $(M,\omega)$ be arbitrary compact symplectic manifold without
boundary. Suppose that $G$ is an autonomous Hamiltonian such that
\roster
\item it has no non-constant contractible periodic orbits ``of period
less than one'',
\item it has a maximum and a minimum that are generically under-twisted
\item all of its critical points are non-degenerate in the Floer
theoretic sense (i.e., the linearized flow of $X_G$ at each critical
point has only the zero as a periodic orbit).
\endroster
Then the one parameter group $\phi_G^t$ is length minimizing in
its homotopy class with fixed ends for $ 0 \leq t \leq 1$.
\endproclaim
A similar result with slightly different assumptions and statements
was proven by McDuff-Slimowitz [MS] by a different method around the same time.

There is also a result by Entov [En] for the strongly semi-positive case.
With some additional restriction on the manifold $(M,\omega)$, we can
remove the condition (3) which we will study elsewhere.

As remarked in [MS] before, the apparently weaker condition ``of
period less than one'' than ``of period less than or equal to
one'' does not give rise to a stronger result. This is because
once we have proven the length minimizing property under the
phrase ``of period less than or equal to one'', the improvement
under the former phrase in Theorem I follows by an approximate
argument as in [Lemma 5.1, Oh3].

We call two Hamiltonians $G$ and $F$ equivalent if
there exists a family $\{F^s\}_{0\leq s\leq 1}$ such that
$$
\phi^1_{F^s} = \phi^1_G
$$
for all $s \in [0,1]$. We denote $G \sim F$ in that case and say
that two Hamiltonian paths $\phi^t_G$ and $\phi^t_F$ are homotopic
to each other with fixed ends, or just homotopic to each other
when there is no danger of confusion.

\definition{Definition 1.1} A Hamiltonian $H$ is called {\it
quasi-autonomous} if there exists two points $x_-, \, x_+ \in M$
such that
$$
H(x_-,t) = \min_x H(x,t), \quad H(x_+,t) = \max_x H(x,t)
$$
for all $t\in [0,1]$.
\enddefinition

We now recall Ustilovsky-Lalonde-McDuff's necessary condition on
the stability of geodesics. Ustilovsky [U] and Lalonde-McDuff [LM]
proved that for a generic $\phi$ in the sense that all its fixed
points are isolated, any stable geodesic $\phi_t, \, 0 \leq t \leq
1$ from the identity to $\phi$ must have at least two fixed points
which are under-twisted.

\medskip

\definition{\bf Definition 1.2} Let $H: M \times [0,1] \to \R$ be a
Hamiltonian which is not necessarily time-periodic and $\phi_H^t$
be its Hamiltonian flow. \par

\roster \item We call a point $p\in M$ a {\it time $T$ periodic
point} if $\phi_H^T(p)=p$. We call $t \in [0,T] \mapsto
\phi_H^t(p)$ a {\it contractible time $T$-periodic orbit} if it is
contractible. \par

\item When $H$ has a fixed critical point $p$ over $t \in
[0,T]$, we call $p$ {\it over-twisted} as a time $T$-periodic
orbit if its linearized flow $d\phi_H^t(p); \, t\in [0,T]$ on
$T_pM$ has a closed trajectory of period less than or equal to
$T$. Otherwise we call it {\it under-twisted}.
If in addition the linearized flow has only the origin as the
fixed point, then we call the fixed point {\it generically under-twisted}.
\endroster
\enddefinition
Here we follow the terminology used in [KL] for the ``generically
under-twisted''.
Note that under this definition of the under-twistedness,
under-twistedness is $C^2$-stable property of the Hamiltonian $H$.

The following is the main result of the present paper, which
improves Theorem I by replacing the phrase ``of period less than
(or equal to) one'' by ``of period one''. This is motivated by a
recent result [KL] of Kerman and Lalonde who first studied the
length minimizing property of the Hamiltonian paths under the
phrase ``of period one'' instead of ``of period less than (or
equal to) one'' on the symplectically aspherical case, with the
same kind of chain level Floer theory as in [Oh3], but specialized
to the {\it symplectically aspherical} case. In this case, the
condition (3) is not needed and the phrase ``in its homotopy
class'' can be replaced by ``among all paths'' as proved in [KL].
We refer readers to [KL] for some explanation on the significance
of such improvement.

\proclaim{Theorem II} Suppose that $G$ is an autonomous
Hamiltonian as in Theorem I except the condition (1) is replaced by
\smallskip

$(1^\prime)\quad$ it has no non-constant contractible periodic orbits ``of
period one''
\smallskip

Then the one parameter group $\phi_G^t$ is length minimizing in
its homotopy class with fixed ends for $ 0 \leq t \leq 1$.
\endproclaim
From now on, we will always assume, unless otherwise said, that
the Hamiltonian functions are normalized so that
$$
\int_M H_t\, d\mu = 0.
$$
When we use a Hamiltonian which is not normalized, we will
explicitly mention it.

Our proof of Theorem II will be again based on the chain level
Floer theory from [Oh3,5], but this time incorporating usage of
the spectral invariants that the author constructed in [Oh5] a
year after the paper [Oh3] appeared. One crucial additional
ingredient in this chain level Floer theory that plays an
important role in our proof of Theorem II is Kerman-Lalonde's
lemma [Proposition 5.2, KL] (see [KL] or \S 4 for detailed account
of this).

In the present paper, in addition to the proof of Theorem II,
using the spectral invariant $\rho(H;1)$ that was constructed in
[Oh5], we provide a much simpler and more elegant scheme than the
one used in [Oh3] for the whole study of length minimizing
property. We note that there has been a general scheme, the so
called {\it energy-capacity inequality}, for the study of length
minimizing property used by Lalonde-McDuff [LM]. Our scheme
belongs to the category of this general scheme. In this respect,
we will state a simple criterion for the length minimizing
property of general quasi-autonomous Hamiltonian paths in terms of
$\rho(\cdot,1)$ on arbitrary closed symplectic manifolds. This
criterion was implicitly used in [Proposition 5.3, Oh3]  without
referring to the spectral invariant. A similar criterion was used
by Hofer [H2] and Bialy-Polterovich [BP] for the compactly
supported Hamiltonians in $\R^{2n}$. Bialy and Polterovich also
predicted existence of similar criterion in general [Remark 1.5,
BP]. The present paper confirms their prediction on arbitrary
closed symplectic manifolds by using the {\it selector}
$\rho(\cdot ;1)$ in their terminology.

To describe this criterion, we rewrite the Hofer norm into
$$
\|H\| = E^-(H) + E^+(H)
$$
where $E^\pm$ are the negative and positive parts of the Hofer
norms defined by
$$
\align
E^-(H) & = \int_0^1 -\min H\, dt \\
E^+(H) & = \int_0^1 \max H\, dt.
\endalign
$$
These are called the {\it negative Hofer-length} and the {\it
positive Hofer-length} of $H$ respectively. We will consider them
separately as usual by now. First note
$$
E^+(H) = E^-(\overline H) \tag 1.3
$$
where $\overline H$ is the Hamiltonian generating
$(\phi_H^t)^{-1}$ defined by
$$
\overline H(t,x) = - H(t, \phi_H^t(x)).
$$
Therefore we will focus only on the semi-norm $E^-$.

\proclaim{Theorem III} Let $G: [0,1] \times M \to \R$ be any
quasi-autonomous Hamiltonian that satisfies
$$
\rho(G;1) = E^-(G) \tag 1.4
$$
Then $G$ is negative Hofer-length minimizing in its homotopy class
with fixed ends.
\endproclaim
The proof will be based on the general property of $\rho(\cdot;1)$
that were proved in [Oh5] which we will recall in \S 3. With this
criterion in mind, Theorem II will follow from the homological
essentialness of the two critical values of $\AA_G$
$$
\align
E^-(G) & := \int_0^1 - \min G\, dt \\
E^+(G) & := \int_0^1 \max G\, dt \\
\endalign
$$
for autonomous Hamiltonian paths of the type as in Theorem II.

\proclaim{Theorem IV} Let $G$ be as in Theorem II. Then (1.4)
holds, i.e., we have
$$
\rho(G;1) = E^-(G).
$$
In particular the critical value $E^-(G)$ is homologically
essential in the Floer theoretic sense. The same holds for
$\overline G$.
\endproclaim
The proof of this theorem is an adaptation of the proof of
Proposition 7.11 (Non-pushing down lemma II) [Oh3] to the current
setting. We will clarify the role of spectral invariant
$\rho(G;1)$ here in its proof.

Finally we would like to compare the scheme of [Oh3] and the
scheme used in the present paper. Both schemes are based on the
mini-max theory via the chain level Floer theory. However while we
explicitly use the chain level Floer theory, more specifically use
sophisticated moving-around of the Floer semi-infinite cycles via
delicate choice of homotopies in [Oh3], these are mostly hidden in
the present paper. This is because we have written the paper [Oh5]
after [Oh3] which provides construction of spectral invariants
whose general properties already reflect this chain level Floer
theory. Furthermore by doing so, we have greatly simplified and
clarified the schemes that we use in [Oh3]. One should note that
statements of the above theorems do not explicitly involve the
Floer theory at all. For example, the Hamiltonian $G$ in Theorem
III is not required to be time one-periodic (see the end of \S 3).
But the Floer theory is implicit and subsumed in the definition of
the spectral invariant $\rho(\cdot; 1)$ in [Oh5]. This may open up
a possibility of suppressing a large part of analytic arguments of
the Floer theory in its application to the study of Hofer's
geodesics or more generally of the Hamiltonian diffeomorphism
group. We will investigate further applications of spectral
invariants to other problems related to the Hamiltonian
diffeomorphism group elsewhere.
\medskip

We would like to thank F. Lalonde for some useful comments on the
previous version of the present paper which lead to clarification
of bibliographical citations, especially concerning the
energy-capacity inequality.

\head \bf \S 2. Preliminaries
\endhead

Let $\Omega_0(M)$ be the set of contractible loops and
$\widetilde\Omega_0(M)$ be its standard covering space in the
Floer theory.  Note that the universal covering space of
$\Omega_0(M)$ can be described as the set of equivalence classes
of the pair $(\gamma, w)$ where $\gamma \in \Omega_0(M)$ and $w$
is a map from the unit disc $D=D^2$ to $M$ such that $w|_{\partial
D} = \gamma$: the equivalence relation to be used is that
$[\overline w \# w^\prime]$ is zero in $\pi_2(M)$. We say that
$(\gamma,w)$ is {\it $\Gamma$-equivalent} to $(\gamma,w^\prime)$
iff
$$
\omega([w'\# \overline w]) = 0 \quad \text{and }\, c_1([w\#
\overline w]) = 0 \tag 2.1
$$
where $\overline w$ is the map with opposite orientation on the
domain and $w'\# \overline w$ is the obvious glued sphere. And
$c_1$ denotes the first Chern class of $(M,\omega)$. We denote by
$[\gamma,w]$ the $\Gamma$-equivalence class of $(\gamma,w)$, by
$\widetilde\Omega_0(M)$ the set of $\Gamma$-equivalence classes
and by $\pi: \widetilde \Omega_0(M) \to \Omega_0(M)$ the canonical
projection. We also call $\widetilde \Omega_0(M)$ the
$\Gamma$-covering space of $\Omega_0(M)$. The unperturbed action
functional $\AA_0: \widetilde \Omega_0(M) \to \R$ is defined by
$$
\AA_0([\gamma,w]) = -\int w^*\omega. \tag 2.2
$$
Two $\Gamma$-equivalent pairs $(\gamma,w)$ and $(\gamma,w^\prime)$
have the same action and so the action is well-defined on
$\widetilde\Omega_0(M)$. When a periodic Hamiltonian $H:M \times
(\R/\Z) \to \R$ is given, we consider the functional $\AA_H:
\widetilde \Omega(M) \to \R$ defined by
$$
\AA_H([\gamma,w])= -\int w^*\omega - \int H(\gamma(t),t)dt
$$
We would like to note that {\it under this convention the maximum
and minimum are reversed when we compare the action functional
$\AA_G$ and the (quasi-autonomous) Hamiltonian $G$}.
One should compare our convention with those used in [Po] or [KL]
where they use the action functional defined by
$$
\AA_H([\gamma,w])= -\int w^*\omega + \int H(\gamma(t),t)\, dt.
$$
Together with their change of the sign in the definition of the
Hamiltonian vector field $X_H$
$$
\dot \iota_{X_H}\omega = -dH,
$$
the difference between the two conventions will be cancelled away
if one makes the substitution of the Hamiltonian
$$
H \longleftrightarrow \widetilde H: \quad \widetilde H(t,x) := -
H(1-t, x).
$$

We denote by $\text{Per}(H)$ the set of periodic orbits of $X_H$.
\medskip

\n{\bf Definition 2.1 \, [Action Spectrum].}  We define the {\it
action spectrum} of $H$, denoted as $\hbox{\rm Spec}(H) \subset
\R$, by
$$
\hbox{\rm Spec}(H) := \{\AA_H(z,w)\in \R ~|~ [z,w] \in
\widetilde\Omega_0(M), z\in \text {Per}(H) \},
$$
i.e., the set of critical values of $\AA_H: \widetilde\Omega(M)
\to \R$. For each given $z \in \text {Per}(H)$, we denote
$$
\hbox{\rm Spec}(H;z) = \{\AA_H(z,w)\in \R ~|~ (z,w) \in
\pi^{-1}(z) \}.
$$

Note that $\text {Spec}(H;z)$ is a principal homogeneous space
modelled by the period group of $(M,\omega)$
$$
\Gamma_\omega = \Gamma(M,\omega) := \{ \omega(A)~|~ A \in \pi_2(M)
\}
$$
and
$$
\hbox{\rm Spec}(H)= \cup_{z \in \text {Per}(H)}\text {Spec}(H;z).
$$
Recall that $\Gamma_\omega$ is either a discrete or a countable
dense subset of $\R$. It is trivial, i.e., $\Gamma_\omega = \{0\}$
in the weakely exact case. The following lemma was proved in
[Oh3].

\proclaim\nofrills{Lemma 2.2. [Lemma 2.2, Oh3]}~ $\hbox{\rm
Spec}(H)$ is a measure zero subset of $\R$.
\endproclaim

For given $\phi \in {\cal  H }am(M,\omega)$, we denote by $H
\mapsto \phi$ if $\phi^1_H = \phi$, and denote
$$
\HH(\phi) = \{ H ~|~ H \mapsto \phi \}, \quad \HH_m(\phi) = \{H
\in \HH(\phi) \mid H \, \text{mean normalized} \}.
$$
We say that two Hamiltonians $H$ and $K$ are equivalent if they
are connected by one parameter family of Hamiltonians
$\{F^s\}_{0\leq s\leq 1}$ such that $F^s \mapsto \phi$ i.e.,
$$
\phi_{F^s}^1 = \phi \tag 2.3
$$
for all $s
\in [0,1]$. We denote by $\widetilde \phi = [\phi, H] = [H]$
the equivalence class of $H$. Then
the universal covering space $\widetilde{{\cal  H }am}(M,\omega)$
of ${\cal  H }am(M,\omega)$ is realized by the set of such
equivalence classes.

Let $F,G \mapsto \phi$ and denote
$$
f_t = \phi_F^t, \, g_t = \phi_G^t,\, \text{and }\,  h_t = f_t
\circ g_t^{-1}.
$$
Note that $h= \{h_t\}$ defines a loop in $\HH am(M,\omega)$ based
at the identity. Suppose $F\sim G$ so there exists a family
$\{F^s\}_{0\leq s \leq 1} \subset \HH(\phi)$ with $F_1 =F$ and
$F_0 = G$ that satisfies (2.3). In particular $h$ defines a
contractible loop.

The following is proved in [Oh4] (see [Sc] for
the symplectically aspherical case where the action functional is
single-valued. In this case Schwarz [Sc] proved that the
normalization works on $\HH am(M,\omega)$ not just on
$\widetilde{\HH am} (M,\omega)$ as long as $F, \, G \mapsto \phi$,
without assuming $F\sim G$).

\proclaim{Proposition 2.3 [Theorem I, Oh4]} Let $F,\, G \in
\HH_m(\phi)$ with $F\sim G$. Then we have
$$
\text{\rm Spec}(G) = \text{\rm Spec}(F)
$$
as a subset of $\R$.
\endproclaim

\head\bf \S 3. Chain level Floer theory and spectral invariants
\endhead

In this section, we will briefly recall the basic chain level
operators in the Floer theory, and the definition and basic
properties of $\rho(\cdot,1)$ from [Oh5].

For each given generic time-periodic $H: M \times S^1 \to \R $, we
consider the free $\Q$ vector space over
$$
\text{Crit}\AA_H = \{[z,w]\in \widetilde\Omega_0(M) ~|~ z \in
\text{Per}(H)\}. \tag 3.1
$$
To be able to define the Floer boundary operator correctly, we
need to complete this vector space downward with respect to the
real filtration provided by the action $\AA_H([z,w])$ of the
element $[z, w]$ of (3.1). More precisely, following [Oh3], we
introduce
\medskip

\definition{Definition 3.1} (1) We call the formal sum
$$
\beta = \sum _{[z, w] \in \text{Crit}\AA_H} a_{[z, w]} [z,w], \,
a_{[z,w]} \in \Q \tag 3.2
$$
a {\it Floer Novikov chain} if there are only finitely many
non-zero terms in the expression (3.2) above any given level of
the action. We denote by $CF(H)$ the set of Novikov chains. We
often simply call them {\it Floer chains}, especially when we do
not need to work on the covering space $\widetilde\Omega_0(M)$ as
in the weakly exact case.

(2) Two Floer chains $\alpha$ and $\alpha'$ are said to be {\it
homologous} to each other if they satisfy
$$
\alpha' = \alpha + \part \gamma
$$
for some Floer chain $\gamma$. We call $\beta$ a {\it Floer cycle}
if $\part \beta = 0$.

(3) Let $\beta$ be a Floer chain in $CF(H)$. We define and denote
the {\it level} of the chain $\beta$ by
$$
\lambda_H(\beta) =\max_{[z,w]} \{\AA_H([z,w] ~|~a_{[z,w]}  \neq
0\, \text{in }\, (3.2) \} \tag 3.3
$$
if $\beta \neq 0$, and just put $\lambda_H(0) = +\infty$ as usual.

(4) We say that $[z,w]$ is a {\it generator} of or {\it
contributes} to $\beta$ and denote
$$
[z,w] \in \beta
$$
if $a_{[z,w]} \neq 0$.
\enddefinition

Let $J = \{J_t\}_{0\leq t \leq 1}$ be a periodic family of
compatible almost complex structures on $(M,\omega)$.

For each given such periodic pair $(J, H)$, we define the boundary
operator
$$
\part: CF(H) \to CF(H)
$$
considering the perturbed Cauchy-Riemann equation
$$
\cases
\frac{\part u}{\part \tau} + J\Big(\frac{\part u}{\part t}
- X_H(u)\Big) = 0\\
\lim_{\tau \to -\infty}u(\tau) = z^-,  \lim_{\tau \to
\infty}u(\tau) = z^+ \\
\endcases
\tag 3.4
$$
This equation, when lifted to $\widetilde \Omega_0(M)$, defines
nothing but the {\it negative} gradient flow of $\AA_H$ with
respect to the $L^2$-metric on $\widetilde \Omega_0(M)$ induced by
the metrics $g_{J_t}: = \omega(\cdot, J_t\cdot)$ . For each given
$[z^-,w^-]$ and $[z^+,w^+]$, we define the moduli space
$$
\MM_{(J,H)}([z^-,w^-],[z^+,w^+])
$$
of solutions $u$ of (3.3) satisfying
$$
w^-\# u \sim w^+. \tag 3.5
$$
$\part$ has degree $-1$ and satisfies $\part\circ \part = 0$.

When we are given a family $(j, \HH)$ with $\HH = \{H^s\}_{0\leq s
\leq 1}$ and $j = \{J^s\}_{0\leq s \leq 1}$, the chain
homomorphism
$$
h_{(j,\HH)}: CF(H^0) \to CF(H^1)
$$
is defined by the non-autonomous equation
$$
\cases \frac{\part u}{\part \tau} +
J^{\rho_1(\tau)}\Big(\frac{\part u}{\part t}
- X_{H^{\rho_2(\tau)}}(u)\Big) = 0\\
\lim_{\tau \to -\infty}u(\tau) = z^-,  \lim_{\tau \to
\infty}u(\tau) = z^+
\endcases
\tag 3.6
$$
where $\rho_i, \, i= 1,2$ is functions of the type $\rho :\R \to
[0,1]$,
$$
\align
\rho(\tau) & = \cases 0 \, \quad \text {for $\tau \leq -R$}\\
                    1 \, \quad \text {for $\tau \geq R$}
                    \endcases \\
\rho^\prime(\tau) & \geq 0
\endalign
$$
for some $R > 0$. We denote by
$$
\MM^{(j,\HH)}([z^-,w^-],[z^+,w^+])
$$
or sometimes with $j$ suppressed the set of solutions of (3.6)
that satisfy (3.5). The chain map $h_{(j,\HH)}$ is defined
similarly as $\part$ using this moduli space instead.
$h_{(j,\HH)}$ has degree 0 and satisfies
$$
\part_{(J^1,H^1)} \circ h_{(j,\HH)} = h_{(j,\HH)} \circ
\part_{(J^0,H^0)}.
$$

Finally,  when we are given a homotopy $(\overline j, \overline
\HH)$ of homotopies with $\overline j = \{j_\kappa\}_{0 \leq
\kappa \leq 1}$, $\overline\HH = \{\HH_\kappa\}_{0\leq \kappa \leq
1}$, consideration of the parameterized version of (3.5) for $ 0
\leq \kappa \leq 1$ defines the chain homotopy map
$$
\widetilde H : CF(H^0) \to CF(H^1)
$$
which has degree $+1$ and satisfies
$$
h_{(j_1, \HH_1)} - h_{(j_0,\HH_0)} = \part_{(J^1,H^1)} \circ
\widetilde H + \widetilde H \circ \part_{(J^0,H^0)}.
$$
Although we will not use this operator explicitly in the present paper,
we have recalled them just for completeness' sake.

The following lemma has played a fundamental role in [Ch],
[Oh1-3,5] and by now become well-known among the experts and can
be proven by a straightforward calculation (see e.g., [Proposition
3.2, Oh3] for its proof).

\proclaim{Lemma 3.2} Let $H, K$ be any Hamiltonian not necessarily
non-degenerate and $j = \{J^s\}_{s \in [0,1]}$ be any given
homotopy and $\HH^{lin} = \{H^s\}_{0\leq s\leq 1}$ be the linear
homotopy $H^s = (1-s)H + sK$. Suppose that (3.5) has a solution
satisfying (3.6). Then we have the identity
$$
\align \AA_F([z^+,w^+]) & - \AA_H([z^-,w^-]) \\
& = - \int \Big|\dudtau \Big|_{J^{\rho_1(\tau)}}^2 -
\int_{-\infty}^\infty \rho_2'(\tau)\Big(F(t,u(\tau,t)) -
H(t,u(\tau,t))\Big) \, dt\,d\tau  \tag 3.7
\endalign
$$
\endproclaim

Now we recall the definition and some basic properties of spectral
invariant $\rho(H;a)$ from [Oh5]. We refer readers to [Oh5] for
the complete discussion on general properties of $\rho(H;a)$.

\proclaim{Definition \& Theorem 3.3 [Oh5]} For any given quantum
cohomology class $0 \neq a \in QH^*(M)$, we have a continuous
function denoted by
$$
\rho_a=\rho(\cdot; a): C_m^0([0,1] \times M) \to \R
$$
such that for two $C^1$ functions $H \sim K$ we have
$$
\rho(H;a) = \rho(K;a) \tag 3.8
$$
for all $a \in QH^*(M)$.
Let $\widetilde \phi, \, \widetilde \psi \in
\widetilde{Ham}(M,\omega)$ and $a \neq 0 \in QH^*(M)$. We define
the map
$$
\rho: \widetilde{Ham}(M,\omega) \times QH^*(M) \to \R
$$
by $\rho(\widetilde\phi;a): = \rho(H;a)$.
\endproclaim

Now we focus on the invariant $\rho(\widetilde \phi; 1)$ for $1
\in QH^*(M)$. We first recall the following quantities
$$
\align E^-(\widetilde \phi) & = \inf_{[\phi,H] = \widetilde \phi}
E^-(H) \tag 3.9\\
E^+(\widetilde \phi) & = \inf_{[\phi,H] = \widetilde \phi} E^+(H).
\tag 3.10
\endalign
$$
The quantities
$$
\rho^{\pm}(\phi) := \inf_{\pi(\widetilde \phi) = \phi}
E^\pm(\widetilde \phi)
$$
then define pseudo-norms on $\HH am(M,\omega)$. It is still an
open question whether $\rho^\pm$ are non-degenerate.

\proclaim{Proposition 3.4 [Theorem II, Oh5]} Let $(M,\omega)$ be
arbitrary closed symplectic manifold. We have
$$
\rho(\widetilde \phi;1) \leq E^-(\widetilde \phi), \quad
\rho(\widetilde \phi^{-1};1) \leq E^+(\widetilde \phi). \tag 3.11
$$
In particular, we have
$$
\rho(H;1) \leq E^-(H), \quad \rho(\overline H; 1) \leq E^+(H) \tag
3.12
$$
for any Hamiltonian $H$.
\endproclaim
For the exact case, the inequality (3.12) had been  earlier proven
in [Oh1,2] in the context of Lagrangian submanifolds and in [Sc]
in for the Hamiltonian diffeomorphim. Now the following theorem
(Theorem III) is an immediate consequence of Theorem 3.3 and
Proposition 3.4.

\proclaim{Theorem 3.5} Let $G: [0,1] \times M \to \R$ be a
quasi-autonomous Hamiltonian. Suppose that $G$ satisfies
$$
\rho(G;1) = E^-(G) \tag 3.13
$$
Then $G$ is negative Hofer-length minimizing in its homotopy class
with fixed ends.
\endproclaim
\demo{Proof} Let $F$ be any Hamiltonian with $F \sim G$. Then we
have a string of equalities and inequality
$$
E^-(G) = \rho(G;1) = \rho(F;1) \leq E^-(F)
$$
from (3.13), (3.8) for $a =1$, (3.12) respectively. This finishes
the proof. \qed\enddemo

So far in this section, we have presumed that the Hamiltonians are
time one-periodic. Now we explain how to dispose the periodicity
and extend the definition of $\rho(H;a)$ for arbitrary time
dependent Hamiltonians $H: [0,1] \times M \to \R$. Note that it is
obvious that the semi-norms $E^\pm(H)$ and $\|H\|$ are defined
without assuming the periodicity. For this purpose, the following
lemma from [Oh3] is important. We leave its proof to readers or to
[Oh3].

\proclaim{Lemma 3.6 [Lemma 5.2, Oh3]} Let $H$ be a given
Hamiltonian $H : [0,1] \times M \to \R$ and $\phi = \phi_H^1$ be
its time-one map. Then we can re-parameterize $\phi_H^t$ in time
so that the re-parameterized Hamiltonian $H'$ satisfies the
following properties: \roster
\item $\phi_{H'}^1 = \phi_H^1$
\item $H' \equiv 0$ near $t = 0, \, 1$ and in particular $H'$ is
time periodic
\item Both $E^\pm(H' - H)$ can be made as small as we want
\item If $H$ is quasi-autonomous, then so is $H'$
\item For the Hamiltonians $H', \, H''$ generating any two such
re-parameterizations of $\phi_H^t$, there is canonical one-one
correspondences between $\text{Per}(H')$ and $\text{Per}(H'')$,
and $\text{Crit }\AA_{H'}$ and $\text{Crit }\AA_{H''}$ with their
actions fixed .

\endroster
Furthermore this re-parameterization is canonical with the
``smallness'' in (3) can be chosen uniformly over $H$ depending
only on the $C^0$-norm of $H$.
\endproclaim

Using this lemma, we can now define $\rho(H;a)$ for arbitrary $H$
by
$$
\rho(H;a): = \rho(H';a)
$$
where $H'$ is the Hamiltonian generating the canonical
re-parameterization of $\phi_H^t$ in time provided in Lemma 3.6.
It follows from (3.8) that this definition is well-defined because
any such re-parameterizations are homotopic to each other with
fixed ends. This being said, we will always assume that our
Hamiltonians are time one-periodic without mentioning further in
the rest of the paper.

\head \bf \S 4. Construction of fundamental Floer cycles
\endhead
In this section and the next, we will prove the following result
(Theorem IV). This in particular proves homologically
essentialness of the critical value
$$
E^-(G) = \int_0^1 - \min G\, dt \tag 4.1
$$
of $\AA_G$.

Note that the hypotheses on $G$ in Theorem IV
already makes it regular in the Floer theory
and so we can define the Floer complex of $G$ {\it without
doing any perturbation on it}.
The proof will use the chain level Floer theory as in [Oh3].

For the proof of Theorem IV, we need to unravel the definition of
$\rho(G;1)$ from [Oh5] in general for arbitrary Hamiltonians $G$.
First for generic (one periodic) Hamiltonians $G$, we consider the
Floer homology class dual to the quantum cohomology class $1 \in
H^*(M) \subset QH^*(M)$, which we denote by $1^\flat$ following
the notation of [Oh5] and call the {\it semi-infinite} fundamental
class of $M$. Then according to [Definition 5.2 \& Theorem 5.5,
Oh5], we have
$$
\rho(G;1) = \inf\{\lambda_G(\gamma) \mid \gamma \in \ker\part_G
\subset CF(G)\, \text{with }\, [\gamma] = 1^\flat \}. \tag 4.2
$$
Then $\rho$ is extended to arbitrary Hamiltonians by continuity in
$C^0$-topology. Therefore to prove (4.1), we need to construct
cycles $\gamma$ with $[\gamma] = 1^\flat$ whose level
$\lambda_G(\gamma)$ become arbitrarily close to $E^-(G)$. In fact,
this was one of the most crucial observations exploited in [Oh3],
without being phrased in terms of the invariant $\rho(G;1)$
because at the time of writing of [Oh3] construction of spectral
invariants in the level of [Oh5] was not carried out yet.

Instead this point was expressed in terms of the existence theorem
of certain Floer's continuity equation over the linear homotopy
(see [Proposition 5.3, Oh3]). Then the author proved the
existence result by proving
homological essentialness of the critical value
$$
E^-(G) = \int_0^1 - \min G\, dt.
$$
The proof relies on a construction of `effective' fundamental
Floer cycles dual to the quantum cohomology class $1$. In [Oh3],
for a suitably chosen Morse function $f$ and for sufficiently
small $\e$, we transferred the fundamental Morse cycle of $\e f$
$$
\alpha_{\e f} := \sum_i^\ell a_{[p_i,w_{p_i}]} [p_i,w_{p_i}] \tag
4.3
$$
to a Floer cycle of $G$ over the {\it adiabatic} homotopy along a
piecewise linear path
$$
\e f \mapsto \e_0 G^{\e_0} \mapsto G \tag 4.4
$$
where $w_p: D^2 \mapsto M$ denote the constant disc $w_p\equiv p$,
and proved the following two facts (see Proposition 7.11 [Oh3]):
\roster
\item the transferred cycle has the level $E^-(G)$ and
\item the cycle cannot be pushed further down under the Cauchy-Riemann
flow
\endroster
under the hypotheses as in Theorem I [Oh3] stated in the
introduction, not just for autonomous but for general
quasi-autonomous Hamiltonians. Now we are ready to introduce the
following fundamental concept of homological essentialness in the
chain level theory, which is already implicitly present in the
series of Non-pushing down lemmas in [Oh3]. As we pointed out in
[Oh3,5], this concept is the heart of the matter in the chain
level theory. In the terminology of [Oh5], the level of any {\it
tight} Floer Novikov cycle of $G$ lies in the {\it essential
spectrum} $\text{spec }G \subset \text{Spec }G$ i.e., realizes the
value $\rho(G;a)$ for some $a \in QH^*(M;\Q)$.

\definition{Definition 4.1} We call a Floer cycle $\alpha \in CF(H)$
{\it tight} if it satisfies the following non-pushing down
property under the Cauchy-Riemann flow (3.4): for any Floer cycle
$\alpha' \in CF(H)$ homologous to $\alpha$ (in the sense of
Definition 3.1 (2)), it satisfies
$$
\lambda_H(\alpha') \geq \lambda_H(\alpha). \tag 4.5
$$
\enddefinition

Now we will attempt to construct a {\it tight} fundamental Floer
cycle of $G$ whose level is precisely $E^-(G)$. As a first step,
we will construct a fundamental cycle of $G$ whose level is
$E^-(G)$ but which may not be tight in general. We choose a Morse
function $f$ such that $f$ has the unique global minimum point
$x^-$ and
$$
f(x^-) = 0, \quad f(x^-) < f(x_j) \tag 4.6
$$
for all other critical points $x_j$. Then we choose a fundamental
Morse cycle
$$
\alpha =\alpha_{\e f} = [x^-,w_{x^-}] + \sum_j a_j [x_j,w_{x_j}]
$$
as in [Oh3] where $x_j \in \text{Crit }_{2n} (-f)$. Recall that
the {\it positive} Morse gradient flow of $\e f$ corresponds to
the {\it negative} gradient flow of $\AA_{\e f}$ in our
convention.

Considering Floer's homotopy map $h_\LL$ over the linear path
$$
\LL: \, s\mapsto (1-s) \e f +s H
$$
for sufficiently small $\e > 0$, we transfer the above fundamental
Morse cycle $\alpha$ and define a fundamental Floer cycle of $H$
by
$$
\alpha_H: = h_{\LL}(\alpha) \in CF(H). \tag 4.7
$$
We call this particular cycle {\it the canonical fundamental Floer
cycle} of $H$. Recently Kerman and Lalonde [KL] proved the
following important property of this fundamental cycle.  Partly
for the reader's convenience and since [KL] only deals with the
aspherical case and our setting is slightly different from [KL],
we give a complete proof here adapting that of [Proposition 5.2,
KL] to our setting of Floer Novikov cycles.

\proclaim{Proposition 4.2 (Compare with [Proposition 5.2, KL])}
Suppose that  $H$ is a generic one-periodic Hamiltonian such that
$H_t$ has the unique  non-degenerate global minimum $x^-$ which is
fixed and under-twisted for all $t \in [0,1]$. Suppose that $f: M
\to \R$ is a Morse function such that $f$ has the unique global
minimum point $x^-$ and $f(x^-)=0$. Then the canonical fundamental
cycle has the expression
$$
\alpha_H = [x^-, w_{x^-}] + \beta \in CF(H)\tag 4.8
$$
for some Floer Novikov chain $\beta \in CF(H)$ with the inequality
$$
\lambda_H(\beta) < \lambda_H([x^-,w_{x^-}]) = \int_0^1 - H(t,
x^-)\, dt. \tag 4.9
$$
In particular its level satisfies
$$
\align \lambda_H(\alpha_H) & = \lambda_H([x^-,w_{x^-}])\tag 4.10\\
& = \int_0^1 - H(t, x^-)\, dt = \int_0^1 -\min H\, dt.
\endalign
$$
\endproclaim

The proof is based on the following simple fact (see the proof of
[Proposition 5.2, KL]). Again we would like to call reader's
attention on the signs due to the different convention we are
using from [KL]. Other than that, we follow the notations from
[KL] in this lemma. To make sure that the different conventions
used in [KL] and here do not cause any problem, we here provide
details of the proof of this lemma.

\proclaim{Lemma 4.3} Let $H$ and $f$ as in Proposition 4.4. Then
for all sufficiently small $\e > 0$, the function $G^H$ defined by
$$
G^H(t,x) = H(t, x^-) + \e f
$$
satisfies
$$
\aligned G^H(t,x^-) & = H(t, x^-) \\
G^H(t,x)& \leq H(t,x)
\endaligned
\tag 4.11
$$
for all $(t,x)$ and equality holds only at $x^-$.
\endproclaim
\demo{Proof} Since $H_t$ has the fixed non-degenerate minimum at
$x^-$ for all $t \in [0,1]$, it follows from a parameterized
version of the Morse lemma that there exists a local coordinates
$(U, y_1, \cdots, y_{2n})$ at $x^-$ such that
$$
H(t,x) = H(t,x^-) + \sum_{i,j = 1}^{2n} a_{ij}(t,x)y_iy_j
$$
with $a_{ij}(t,x^-)$ is a positive definite matrices for each $t
\in [0,1]$ which depend smoothly on $t$. On the coordinate
neighborhood $U$, we have
$$
\align H(t,x) - G^H(t,x) & = H(t,x) - (H(t,x^-) + \e f(x)) \\
& = \sum_{i,j = 1}^{2n} a_{ij}(t,x)y_iy_j - \e f(x). \tag 4.12
\endalign
$$
Since $f$ has the non-degenerate minimum point at $x^-$ and
$f(x^-) = 0$, it follows from (4.12) that for any sufficiently
small $\e > 0$, we have
$$
H(t,x) - G^H(t,x) \geq 0 \quad \text{for all } \, (t,x) \in [0,1]
\times U \tag 4.13
$$
and {\it equality only at $x = x^-$}, if we choose sufficiently
small $U$. On the other hand, since $x^-$ is the {\it unique}
fixed non-degenerate global minimum of $H$, there exists $\delta >
0$ such
$$
H(t,x) - H(t,x^-) \geq \delta
$$
for all $(t,x) \in [0,1] \times (M \setminus U)$.  If we choose
$\e$ so that $\e \max f \leq {1 \over 2}\delta$, we also have
$$
H(t,x) - G^H(t,x) \geq {1 \over 2}\delta \quad \text{ for all }\,
(t,x) \in [0,1] \times (M \setminus U). \tag 4.14
$$
Combining (4.13) and (4.14), we have finished the proof.
\qed\enddemo

\demo{Proof of Proposition 4.2} Since $x^-$ is a under-twisted
fixed minimum of both $H$ and $f$, we have the Conley-Zehnder
index
$$
\mu_H([x^-, w_{x^-}])  = \mu_{\e f}([x^-,w_{x^-}]) ( = -n)
$$
and so the moduli space $\MM^\LL([x^-,w_{x^-}],[x^-,w_{x^-}])$ has
dimension zero. Let $u \in \MM^\LL([x^-,w_{x^-}],[x^-,w_{x^-}])$.

We note that the Floer continuity equation (3.6) for the linear
homotopy
$$
\LL: s \to (1-s) \e f + s H
$$
is unchanged even if we replace the homotopy by the homotopy
$$
\LL': s \to (1-s) G^H + s H.
$$
This is because the added term $H(t,x^-)$ in $G^H$ to $\e f$ does
not depend on $x \in M$ and so
$$
X_{\e f} \equiv X_{G^H}.
$$
Therefore $u$ is also a solution for the continuity equation (3.6)
under the linear homotopy $\LL'$. Using this,  we derive the
identity
$$
\aligned \int \Big|\dudtau \Big|_{J^{\rho_1(\tau)}}^2\, dt\,d\tau
& =
\AA_{G^H}([x^-,w_{x^-}]) - \AA_H([x^-,w_{x^-}]) \\
\quad & - \int_{-\infty}^\infty \rho'(\tau)\Big(H(t,u(\tau,t))\,
dt\, d\tau - G^H(t,u(\tau,t))\Big) \, dt\,d\tau
\endaligned
\tag 4.15
$$
from (3.7). Since we have
$$
\AA_H([x^-,w_{x^-}]) =\AA_{G^H}([x^-,w_{x^-}]) = \int_0^1 - \min
H\, dt \tag 4.16
$$
and $G^H \leq H$, the right hand side of (4.15) is non-positive.
Therefore we derive that $\MM^\LL([x^-,w_{x^-}],[x^-,w_{x^-}])$
consists only of the constant solution $u \equiv x^-$. This in
particular gives rise to the matrix coefficient of $h_\LL$
satisfying
$$
\langle [x^-,w_{x^-}], h_{\LL}([x^-,w_{x^-}])\rangle = \#
(\MM^\LL([x^-,w_{x^-}],[x^-,w_{x^-}])) = 1.
$$
Now consider any other generator of $\alpha_H$
$$
[z,w] \in \alpha_H \quad \text{with }\,  [z,w] \neq [x^-,w_{x^-}].
$$
By the definitions of $h_\LL$ and $\alpha_H$, there is a generator
$[x,w_x] \in \alpha$ such that
$$
\MM^\LL([x,w_x],[z,w]) \neq \emptyset. \tag 4.17
$$
Then for any $u \in \MM^\LL([x,w_x],[z,w])$, we have the identity
from (3.7)
$$
\aligned  \AA_H([z,w]) - & \AA_{G^H}([x,w_x]) = -\int
\Big|\dudtau \Big|_{J^{\rho_1(\tau)}}^2\, dt\, d\tau \\
& \quad - \int_{-\infty}^\infty \rho'(\tau)\Big(H(t,u(\tau,t)) -
G^H(t,u(\tau,t))\Big) \, dt\,d\tau. \endaligned
$$
Since $-\int \Big|\dudtau \Big|_{J^{\rho_1(\tau)}}^2 \leq 0$, and
$G^H \leq H$,  we have
$$
\AA_H([z,w]) \leq \AA_{G^H}([x,w_x])\tag 4.18
$$
with equality holding only when $u$ is stationary. There are two
cases to consider, one for the case of $x = x^-$ and the other for
$x = x_j$ for $x_j \neq x^-$ for $[x_j,w_{x_j}] \in \alpha$.

For the first case, {\it since we assume $[z,w] \neq
[x^-,w_{x^-}]$}, $u$ cannot be constant and so the strict
inequality holds in (4.18), i.e,
$$
\AA_H([z,w]) < \AA_{G^H}([x^-,w_{x^-}]). \tag 4.19
$$
For the second case, we have the inequality
$$
\AA_H([z,w]) \leq \AA_{G^H}([x_j,w_{x_j}])\tag 4.20
$$
for some $x_j \neq x^-$ with $[x_j,w_{x_j}] \in \alpha$. We note
that (4.6) is equivalent to
$$
\AA_{G^H}([x_j,w_{x_j}]) < \AA_{G^H}([x^-,w_{x^-}]).
$$
This together with (4.20) again give rise to (4.19). On the other
hand we also have
$$
\AA_{G^H}([x^-,w_{x^-}]) = \AA_H([x^-,w_{x^-}])
$$
because $G^H(t,x^-) = H(t,x^-)$ from (4.11). Altogether, we have
proved
$$
\AA_H([z,w]) < \AA_H([x^-,w_{x^-}]) = \int_0^1 -H(t,x^-) \, dt
$$
for any $[z,w] \in \alpha_H$ with $[z,w] \neq [x^-,w_{x^-}]$. This
finishes the proof of (4.9). \qed\enddemo

\definition{Remark 4.4} Note that $G^H$ does not necessarily
satisfy the normalization condition. This causes no problem
because the proof of Proposition 4.4 does not require
normalization condition.
\enddefinition

\head \bf\S 5. The case of autonomous Hamiltonians
\endhead

In this section, we will restrict to the case of autonomous
Hamiltonians $G$ as in Theorem II and prove the following theorem.

\proclaim{Theorem 5.1} Suppose that $G$ is autonomous
as in Theorem II.  Then the canonical
fundamental cycle is tight in the sense of Definition 4.3, i.e.,
$\alpha_G$ satisfies non-pushing down property: for any Floer
Novikov cycle $\alpha \in CF(G)$ homologous to $\alpha_G$, we have
$$
\lambda_{G}(\alpha) \geq \lambda_{G}(\alpha_G). \tag 5.1
$$
In particular, we have
$$
\rho(G;1) = \lambda_G(\alpha_G)= \int_0^1 -\min G = E^-(G). \tag
5.2
$$
\endproclaim
\demo{Proof} The proof is an adaptation of the proof of
Proposition 7.11 (Non-pushing down lemma II) [Oh3].
Note that the conditions in Theorem II in particular
impliy that $G$ is regular in the sense of the Floer theory.

Suppose that $\alpha$ is homologous to $\alpha_G$, i.e.,
$$
\alpha = \alpha_G + \part_G (\gamma) \tag 5.3
$$
for some Floer Novikov chain $\gamma \in CF(G)$. When $G$ is
autonomous and $J\equiv J_0$ is $t$-independent,  there is no
non-stationary $t$-independent trajectory of $\AA_{G}$ landing at
$[x^-,w_{x^-}]$ because any such trajectory comes from the
negative Morse gradient flow of $G$ but $x^-$ is the minimum point
of $G$. Therefore any non-stationary Floer trajectory $u$ landing
at $[x^-,w_{x^-}]$ must be $t$-dependent. Because of the
assumption that $G$ has no non-constant contractible periodic
orbits of period one, any critical points of $\AA_G$ has the form
$$
[x,w]\quad \text{with } \, x \in \text{Crit }G.
$$
Let $u$ be a trajectory starting at $[x,w]$, $x \in \text{Crit }G$
with
$$
\mu([x,w]) - \mu([x^-,w_{x^-}]) = 1, \tag 5.4
$$
and denote by $\MM_{(J_0, G)}([x,w],[x^-,w_{x^-}])$ the
corresponding Floer moduli space of connecting trajectories. The
general index formula shows
$$
\mu([x,w]) = \mu([x,w_{x}]) + 2 c_1([w]). \tag 5.5
$$
We consider two cases separately:  the cases of $c_1([w]) = 0$ or
$c_1([w]) \neq 0$. If $c_1([w]) \neq 0$, we derive from (5.4),
(5.5) that $x \neq x^-$. This implies that any such trajectory
must come with (locally) free $S^1$-action, i.e., the moduli space
$$
\widehat{\MM}_{(J_0,G)}([x,w],[x^-,w_{x^-}]) =
\MM_{(J_0,G)}([x,w],[x^-,w_{x^-}])/\R
$$
and its stable map compactification have a locally free
$S^1$-action {\it without fixed points}. Therefore after a
$S^1$-invariant perturbation $\Xi$ via considering the quotient
Kuranishi structure [FOn] on the quotient space
$\widehat{\MM}_{(J_0,G)}([x,w],[x^-,w_{x^-}])/S^1$, the
corresponding perturbed moduli space
$\widehat{\MM}_{(J_0,G)}([x,w],[x^-,w_{x^-}]; \Xi)$ becomes empty
for a $S^1$-equivariant perturbation $\Xi$. This is because the
quotient Kuranishi structure has the virtual dimension -1 by the
assumption (5.4). We refer to [FHS], [FOn] or [LT] for more
explanation on this $S^1$-invariant regularization process. Now
consider the case $c_1([w]) = 0$. First note that (5.4) and (5.5)
imply that $x \neq x^-$. On the other hand, if $x\neq x^-$, the
same argument as above shows that the perturbed moduli space
becomes empty.

It now follows that there is no trajectory of index 1 that land at
$[x^-,w_{x^-}]$ after the $S^1$-invariant regularization.
Therefore $\part_G(\gamma)$ cannot kill the term $[x^-,w_{x^-}]$
in (5.3) away from the cycle
$$
\alpha_G = [x^-, w_{x^-}] + \beta
$$
in (4.9), and hence we have
$$
\lambda_G(\alpha) \geq \lambda_G([x^-,w_{x^-}]) \tag 5.6
$$
by the definition of the level $\lambda_G$.  Combining (4.10) and
(5.6), we have finished the proof (5.1). \qed\enddemo

\head {\bf References}
\endhead

\enddocument